# Sums of Independent Random Variables for Shanker, Akash, Ishita, Pranav, Rani and Ram Awadh Distributions

Afshin Yaghoubi

Faculty of Mathematics and Computer Science, Amirkabir University of Technology, Tehran, Iran;

Email: afshin.y@aut.ac.ir

**Abstract**

In statistics and probability theory, one of the most important statistics is the sums of random variables. After introducing a probability distribution, determining the sum of *n* independent and identically distributed random variables is one of the interesting topics for authors. This paper presents the probability density functions for the sum of *n* independent and identically distributed random variables such as Shanker, Akash, Ishita, Pranav, Rani, and Ram Awadh. In order to determine all aforementioned distributions, the problem-solving methods are applied which is based on the change-of-variables technique. The *m*th moments for them were also accurately calculated. Besides, the reliability and the mean time to failure of a 1 out of *n* cold standby spare system has also been evaluated under the Lindley components failure time.

**Keywords:** Akash distribution, Ishita distribution, Pranav distribution, Ram Awadh distribution, Rani distribution, Shanker distribution.

1. Introduction

The items such as various wireless communications (e.g., satellite communications, equal-gain receivers, and radars) and insurance are widely known as applications for the sums of random variables (Karagiannidis et al. (2005), Khuong and Kong (2006), Nadarajah (2008), Ramsay (2008) and Nadarajah et al. (2012)). Another application of the sum of random variables is in the field of engineering. For instance, the sum of excess water output in a dam is equal to $X_1 + \cdots + X_n$, where $X_i$ represents the *i*th additional current at position *i* (Moschopoulos (1985)). Moreover, in the field of systems reliability is also used. For example, the lifespan of 1 out of *n* cold standby spare systems is equal to the total lifetime of system components. In other word, if $S_n$ represents the lifetime of the whole system and $X_i$ indicates the time to failure of the *i*th component, then $S_n = X_1 + \cdots X_n$ (C. Wang et al. (2021)).

Recently, a collection of the one-parameter continuous probability density functions (PDFs) has been introduced by Shukla and Shanker (2019). Similar to the Lindley distribution, these functions are also made up of a mixture of exponential distribution with scale parameter $\theta$, $(\exp(\theta))$ and gamma distribution with different shape parameters $n$, and scale parameters $\theta$, $(\text{gamma}(n, \theta))$.

Let consider $f_E(x)$ and $f_G(x)$ as the exponential and the gamma PDFs, respectively. In this case, their combination with a specified ratio of $p$, namely $f(x) = pf_E(x) + (1-p)f_G(x)$ represents the proposed PDFs. The information required to create $f(x)$'s are summarized in Table 1.

**Table 1.** The PDFs of Shanker, Akash, Ishita, Pranav, Rani, and Ram Awadh

| Distributions | $f_E(x)$ | $f_G(x)$ | $p$ | PDFs | By |
|---|---|---|---|---|---|
| Shanker | $\theta e^{-\theta x}$ | $\dfrac{\theta^2 x e^{-\theta x}}{}$ | $\dfrac{\theta^2}{\theta^2+1}$ | $f(x) = \dfrac{\theta^2}{\theta^2+1}(\theta+x)e^{-\theta x}$ | Shanker, 2015 a |
| Akash | $\theta e^{-\theta x}$ | $\dfrac{\theta^3 x^2 e^{-\theta x}}{2!}$ | $\dfrac{\theta^3}{\theta^2+2}$ | $f(x) = \dfrac{\theta^3}{\theta^2+2}(1+x^2)e^{-\theta x}$ | Shanker, 2015 b |
| Ishita | $\theta e^{-\theta x}$ | $\dfrac{\theta^3 x^2 e^{-\theta x}}{2!}$ | $\dfrac{\theta^3}{\theta^3+2}$ | $f(x) = \dfrac{\theta^3}{\theta^3+2}(\theta+x^2)e^{-\theta x}$ | Shanker and Shukla, 2017 |
| Pranav | $\theta e^{-\theta x}$ | $\dfrac{\theta^4 x^3 e^{-\theta x}}{3!}$ | $\dfrac{\theta^4}{\theta^4+6}$ | $f(x) = \dfrac{\theta^4}{\theta^4+6}(\theta+x^3)e^{-\theta x}$ | Shukla, 2018 |
| Rani | $\theta e^{-\theta x}$ | $\dfrac{\theta^5 x^4 e^{-\theta x}}{4!}$ | $\dfrac{\theta^5}{\theta^5+24}$ | $f(x) = \dfrac{\theta^5}{\theta^5+24}(\theta+x^4)e^{-\theta x}$ | Shanker, 2017 |
| Ram Awadh | $\theta e^{-\theta x}$ | $\dfrac{\theta^6 x^5 e^{-\theta x}}{5!}$ | $\dfrac{\theta^6}{\theta^6+120}$ | $f(x) = \dfrac{\theta^6}{\theta^6+120}(\theta+x^5)e^{-\theta x}$ | Shukla, 2018 |

The authors have compared available distributions in Table 1 with the Lindley distribution and the results showed that the performance of the them is better than the Lindley distribution.

The Lindley distribution was introduced by Lindley (1958). This distribution is plentifully used in the reliability theory. Ghitany et al. (2008) examined the Lindley distribution in detail and showed that this distribution in many ways has more flexibility than exponential distribution.

The common methods for obtaining the sums of $n$ IID random variables are the moment generating functions (MGF), characteristic function, Laplace transforms, convolution method, and change-of-variables technique. These methods can be used to find the sum of the distributions which mathematically have a simple PDF (e.g., exponential, gamma, Lindley, and Pareto). For example, Khuong & Kong (2006) by using characteristic function obtained the distribution of sums for the exponential distribution. Ramsay (2008), using the Laplace transform, presented an exact formula for the sum of $n$ independent Pareto random variables. The sum of *n* independent variables that follow the Lindley distribution are obtained by authors such as Zakerzadeh and Dolati (2009), Al-Mutairi et al. (2013) and Hassan (2014). Zakerzadeh and Dolati (2009) provided the distribution of sum for the generalized Lindley (GL) by using the MGF. Al-Mutairi et al. (2013) also investigated the sums of *n* IID Lindley distribution through MGF. Hassan (2014) also obtained it using the convolution integrals.

Distribution of sums for other density functions, like the Lognormal, Weibull, Rayleigh, Burr and so on that have a more complex mathematical form, methods such as of the generalized hypergeometric functions such as the Fox's H-function and the Meijer's G-function are usually used. For example, can be referred to the works Nadarajah et al. (2012), El Bouanani & da Costa (2018), Karagiannidis ae al. (2005) and Beaulieu & Rajwani (2004).

Introducers of the presented PDFs in Table 1 have discussed some features such as descriptive statistics (mean, variance, coefficient of variation, skewness and kurtosis), reliability theory (hazard function, mean residual life and stress-strength reliability models), stochastic ordering, mean deviations, Lorenz curves, entropies and simulation study. But another important aspect, i.e., the sum of *n* independent and identically distributed (IID) random variables for introduced functions, has not been studied. Therefore, here we calculate the distribution of sums for the aforementioned PDFs.

2. **The distribution of sums**

Suppose $X_1, X_2, \ldots, X_n$ denote *n* IID random variables that according to the PDFs of Table 1. If $S_n = X_1 + X_2 + \cdots + X_n$ be the distribution of sums for them, then applying the change-of-variables technique one can obtain the density functions of $S_n$ for Shanker, Akash, Ishita, Pranav, Rani, and Ram Awadh. The results are briefly shown in Table 2.

**Table 2.** Distribution of sums for PDFs of Table 1

| Distributions | PDFs |
|---|---|
| Shanker | $f_{S_n}(x) = \left(\dfrac{\theta^2}{\theta^2+1}\right)^n e^{-\theta x} \sum_{r=0}^{n} \left(\binom{n}{r} \theta^{n-r} \dfrac{x^{n+r-1}}{(n+r-1)!}\right)$ |
| Akash | $f_{S_n}(x) = \left(\dfrac{\theta^3}{\theta^2+2}\right)^n e^{-\theta x} \sum_{r=0}^{n} \left(\binom{n}{r} 2^r \dfrac{x^{n+2r-1}}{(n+2r-1)!}\right)$ |
| Ishita | $f_{S_n}(x) = \left(\dfrac{\theta^3}{\theta^3+2}\right)^n e^{-\theta x} \sum_{r=0}^{n} \left(\binom{n}{r} \theta^{n-r} 2^r \dfrac{x^{n+2r-1}}{(n+2r-1)!}\right)$ |
| Pranav | $f_{S_n}(x) = \left(\frac{\theta^4}{\theta^4 + 6}\right)^n e^{-\theta x} \sum_{r=0}^{n} \binom{n}{r} \theta^{n-r} 6^r \frac{x^{n+3r-1}}{(n+3r-1)!}$ |
| Rani | $f_{S_n}(x) = \left(\frac{\theta^5}{\theta^5 + 24}\right)^n e^{-\theta x} \sum_{r=0}^{n} \binom{n}{r} \theta^{n-r} 24^r \frac{x^{n+4r-1}}{(n+4r-1)!}$ |
| Ram Awadh | $f_{S_n}(x) = \left(\frac{\theta^6}{\theta^6 + 120}\right)^n e^{-\theta x} \sum_{r=0}^{n} \binom{n}{r} \theta^{n-r} 120^r \frac{x^{n+5r-1}}{(n+5r-1)!}$ |

To prove the PDFs of Table 2, assuming $X_i = U_i$, for $i = 1, 2, \ldots, n-1$, and $X_n = S_n - \sum_{i=1}^{n-1} U_i$, can be written joint density of $(U_1, U_2, \ldots, S_n)$ as follows

$$f_{U_1, U_2, \ldots, S_n}(u_1, u_2, \ldots, s_n) = f_{X_1, X_2, \ldots, X_n}\left(u_1, u_2, \ldots, s_n - \sum_{i=1}^{n-1} u_i\right) \times |J|. \qquad (1)$$

where, $|J|$ is the absolute value of the Jacobian.

As the relationship among $X_i$ and $U_i$ is linear, so $|J| = \left|\frac{\partial(X_1, X_2, \ldots, X_n)}{\partial(U_1, U_2, \ldots, S_n)}\right| = 1$. Furthermore, since $X_i$ is independent, hence Eq. (1) for $u_i > 0, s_n > \sum_{i=1}^{n-1} u_i$ is given by

$$f_{U_1, U_2, \ldots, S_n}(u_1, u_2, \ldots, s_n) = f_{X_1}(u_1) f_{X_2}(u_2) \ldots f_{X_n}\left(s_n - \sum_{i=1}^{n-1} u_i\right). \qquad (2)$$

and otherwise, $f_{U_1, U_2, \ldots, S_n}(u_1, u_2, \ldots, s_n) = 0$.

Now, the marginal PDF for $S_n$ for $u_i > 0, s_n > \sum_{i=1}^{n-1} u_i$ is equal to

$$f_{S_n}(x) = \int_0^x \int_0^{x-u_1} \ldots \int_0^{x-\sum_{i=1}^{n-2} u_i} \left(\prod_{i=1}^{n-1} f_{X_i}(u_i)\right) f_{X_n}\left(s_n - \sum_{i=1}^{n-1} u_i\right) du_{n-1} \ldots du_2 du_1. \qquad (3)$$

In Eq. (3), as can be observed, the convolution integrals make the calculations somewhat difficult. Therefore, inductive reasoning is used to solve this problem.

### 2.1. Shanker distribution

By taking the shanker distribution in Eq. (3) and for $x > 0$, the functions of $f_{S_n}(x)$, for $n=2$ and 3 are given by

$$f_{S_2}(x) = \left(\frac{\theta^3}{\theta^2 + 2}\right)^2 e^{-\theta x} \int_0^x (\theta + u_1^2)(\theta + (x - u_1)^2) du_1 = \left(\frac{\theta^3}{\theta^3 + 2}\right)^2 e^{-\theta x} A_1,$$

$$f_{S_3}(x) = \left(\frac{\theta^3}{\theta^3 + 2}\right)^3 e^{-\theta x} \int_0^x \int_0^{x-u_1} (\theta + u_1^2)(\theta + u_2^2)(\theta + (x - u_1 - u_2)^2) du_2 du_1$$

$$= \left(\frac{\theta^3}{\theta^3 + 2}\right)^3 e^{-\theta x} A_2.$$

where $A_1$ and $A_2$ are respectively the result of single and double integrals such that

$$A_1 = \theta^2 x + \theta x^2 + \frac{x^3}{6},$$

$$A_2 = \frac{60\theta^3 x^2}{120} + \frac{60\theta^2 x^3}{120} + \frac{15\theta x^4}{120} + \frac{x^5}{120} = \theta^3 \frac{x^2}{2} + \theta^2 \frac{x^3}{2} + \theta \frac{x^4}{8} + \frac{x^5}{120}.$$

The coefficients of $\frac{x^i}{i!}$ in the phrase of $A_1$ and $A_2$, after some algebraic manipulations, can be rewritten as

$$A_1 = \binom{2}{0}\theta^{2-0}\frac{x^1}{1!} + \binom{2}{1}\theta^{2-1}\frac{x^2}{2!} + \binom{2}{2}\theta^{2-2}\frac{x^3}{3!},$$

$$A_2 = \binom{3}{0}\theta^{3-0}\frac{x^2}{2!} + \binom{3}{1}\theta^{3-1}\frac{x^3}{3!} + \binom{3}{2}\theta^{3-2}\frac{x^4}{4!} + \binom{3}{3}\theta^{3-3}\frac{x^5}{5!}.$$

According to the above results, it can be concluded that $A_{n-1}$ is equals to

$$A_{n-1} = \sum_{r=0}^{n} \binom{n}{r} \theta^{n-r} \frac{x^{n+r-1}}{(n+r-1)!}.$$

In fact, $A_{n-1}$ is the calculation of the ($n$-1)-fold convolution integral, that is

$$\int_0^x \int_0^{x-u_1} \cdots \int_0^{x-\sum_{i=1}^{n-2} u_i} \left( \prod_{i=1}^{n-1}(\theta + u_i) \right)\left( \theta + x - \sum_{i=1}^{n-1} u_i \right) du_{n-1} \ldots du_2 du_1 = A_{n-1}.$$

So, $f_{S_n}(x)$ with Shanker distribution is equals to

$$f_{S_n}(x) = \left(\frac{\theta^2}{\theta^2 + 1}\right)^n e^{-\theta x} \sum_{r=0}^{n} \left(\binom{n}{r} \theta^{n-r} \frac{x^{n+r-1}}{(n+r-1)!}\right). \tag{4}$$

**2.2. Akash distribution**

Replacing the Akash distribution in Eq. (3) and for $x > 0$, the functions of $f_{S_n}(x)$, for $n$=2 and 3 are given by

$$f_{S_2}(x) = \left(\frac{\theta^3}{\theta^2 + 2}\right)^2 e^{-\theta x} \int_0^x (1 + u_1^2)(1 + (x - u_1)^2)du_1 = \left(\frac{\theta^3}{\theta^3 + 2}\right)^2 e^{-\theta x} B_1,$$

$$f_{S_3}(x) = \left(\frac{\theta^3}{\theta^3 + 2}\right)^3 e^{-\theta x} \int_0^x \int_0^{x-u_1} (1 + u_1^2)(1 + u_2^2)(1 + (x - u_1 - u_2)^2) du_2 d du_1$$

$$= \left(\frac{\theta^3}{\theta^3 + 2}\right)^3 e^{-\theta x} B_2.$$

where $B_1$ and $B_2$ are

$$B_1 = \frac{30x}{30} + \frac{20x^3}{30} + \frac{x^5}{30} = x + 2\frac{x^3}{3} + \frac{x^5}{30},$$

$$B_2 = \frac{2520x^2}{5040} + \frac{1260x^4}{5040} + \frac{84x^6}{5040} + \frac{x^8}{5040} = \frac{x^2}{2} + \frac{x^4}{4} + \frac{x^6}{60} + \frac{x^8}{5040}.$$

The $B_1$ and $B_2$ can be rewritten with algebraic calculations as

$$B_1 = \binom{2}{0} 2^0 \frac{x^1}{1!} + \binom{2}{1} 2^1 \frac{x^2}{2!} + \binom{2}{2} 2^2 \frac{x^3}{3!},$$

$$B_2 = \binom{3}{0} 2^0 \frac{x^2}{2!} + \binom{3}{1} 2^1 \frac{x^3}{3!} + \binom{3}{2} 2^2 \frac{x^4}{4!} + \binom{3}{3} 2^3 \frac{x^5}{5!}.$$

According to the above results, it can be concluded that $B_{n-1}$ is equals to

$$B_{n-1} = \sum_{r=0}^{n} \binom{n}{r} 2^r \frac{x^{n+2r-1}}{(n+2r-1)!}.$$

Thus

$$\int_0^x \int_0^{x-u_1} \cdots \int_0^{x-\sum_{i=1}^{n-2} u_i} \left( \prod_{i=1}^{n-1} (1+u_i^2) \right) \left( 1+x-\sum_{i=1}^{n-1} u_i^2 \right) du_{n-1} \ldots du_2 du_1 = B_{n-1}.$$

So, $f_{S_n}(x)$ with Akash distribution is equals to

$$f_{S_n}(x) = \left( \frac{\theta^3}{\theta^2+2} \right)^n e^{-\theta x} \sum_{r=0}^{n} \binom{n}{r} 2^r \frac{x^{n+2r-1}}{(n+2r-1)!}. \tag{5}$$

### 2.3. Ishita distribution

With placement the Ishita distribution in Eq. (3) and for $x > 0$, the functions of $f_{S_n}(x)$, for $n=2$ and 3 are given by

$$f_{S_2}(x) = \left( \frac{\theta^3}{\theta^3+2} \right)^2 e^{-\theta x} \int_0^x (\theta + u_1^2)(\theta + (x-u_1)^2) du_1 = \left( \frac{\theta^3}{\theta^3+2} \right)^2 e^{-\theta x} C_1,$$

$$f_{S_3}(x) = \left( \frac{\theta^3}{\theta^3+2} \right)^3 e^{-\theta x} \int_0^x \int_0^{x-u_1} (\theta + u_1^2)(\theta + u_2^2)(\theta + (x-u_1-u_2)^2) du_2 du_1$$

$$= \left( \frac{\theta^3}{\theta^3+2} \right)^3 e^{-\theta x} C_2.$$

where

$$C_1 = \frac{30x}{30} + \frac{20x^3}{30} + \frac{x^5}{30} = x + 2\frac{x^3}{3} + \frac{x^5}{30},$$

$$C_2 = \frac{2520x^2}{5040} + \frac{1260x^4}{5040} + \frac{84x^6}{5040} + \frac{x^8}{5040} = \frac{x^2}{2} + \frac{x^4}{4} + \frac{x^6}{60} + \frac{x^8}{5040}.$$

Subsequently

$$C_1 = \binom{2}{0} 2^0 \theta^{2-0} \frac{x^1}{1!} + \binom{2}{1} 2^1 \theta^{2-1} \frac{x^3}{3!} + \binom{2}{2} 2^2 \theta^{2-2} \frac{x^5}{5!},$$

$$C_2 = \binom{3}{0} 2^0 \theta^{3-0} \frac{x^2}{2!} + \binom{3}{1} 2^1 \theta^{3-1} \frac{x^4}{4!} + \binom{3}{2} 2^2 \theta^{3-2} \frac{x^6}{6!} + \binom{3}{3} 2^3 \theta^{3-3} \frac{x^8}{8!}.$$

According to the above results, it can be concluded that $C_{n-1}$ is equals to

$$C_{n-1} = \sum_{r=0}^{n} \binom{n}{r} 2^r \theta^{n-r} \frac{x^{n+2r-1}}{(n+2r-1)!}.$$

Therefore

$$\int_0^x \int_0^{x-u_1} \cdots \int_0^{x-\sum_{i=1}^{n-2} u_i} \left(\prod_{i=1}^{n-1} (\theta + u_i^2)\right) \left(\theta + x - \sum_{i=1}^{n-1} u_i^2\right) du_{n-1} \dots du_2 du_1 = C_{n-1}.$$

So, $f_{S_n}(x)$ with Ishita distribution is equals to

$$f_{S_n}(x) = \left(\frac{\theta^3}{\theta^3 + 2}\right)^n e^{-\theta x} \sum_{r=0}^{n} \left(\binom{n}{r} 2^r \theta^{n-r} \frac{x^{n+2r-1}}{(n+2r-1)!}\right). \tag{6}$$

**2.4. Pranav distribution**

Replacement the Pranav distribution in Eq. (3) and for $x > 0$, the functions of $f_{S_n}(x)$, for $n=2$ and 3 are given by

$$f_{S_2}(x) = \left(\frac{\theta^4}{\theta^4 + 6}\right)^2 e^{-\theta x} \int_0^x (\theta + u_1^3)(\theta + (x - u_1)^3) du_1 = \left(\frac{\theta^4}{\theta^4 + 6}\right)^2 e^{-\theta x} D_1,$$

$$f_{S_3}(x) = \left(\frac{\theta^4}{\theta^4 + 6}\right)^3 e^{-\theta x} \int_0^x \int_0^{x-u_1} (\theta + u_1^3)(\theta + u_2^3)(\theta + (x - u_1 - u_2)^3) du_2 du_1$$

$$= \left(\frac{\theta^4}{\theta^4 + 6}\right)^3 e^{-\theta x} D_2.$$

where

$$D_1 = \theta^2 x + \frac{\theta x^4}{2} + \frac{x^7}{140},$$

$$D_2 = \frac{\theta^3 x^2}{2} + \frac{3\theta^2 x^5}{20} + \frac{3\theta x^8}{1120} + \frac{x^{11}}{184800}.$$

Subsequently

$$D_1 = \binom{2}{0} 6^0 \theta^{2-0} \frac{x^1}{1!} + \binom{2}{1} 6^1 \theta^{2-1} \frac{x^4}{4!} + \binom{2}{2} 6^2 \theta^{2-2} \frac{x^7}{7!},$$

$$D_2 = \binom{3}{0} 6^0 \theta^{3-0} \frac{x^2}{2!} + \binom{3}{1} 6^1 \theta^{3-1} \frac{x^5}{5!} + \binom{3}{2} 6^2 \theta^{3-2} \frac{x^8}{8!} + \binom{3}{3} 6^3 \theta^{3-3} \frac{x^{11}}{11!}.$$

According to the above results, it can be concluded that $D_{n-1}$ is equals to

$$D_{n-1} = \sum_{r=0}^{n} \binom{n}{r} 6^r \theta^{n-r} \frac{x^{n+3r-1}}{(n+3r-1)!}.$$

Nevertheless

$$\int_0^x \int_0^{x-u_1} \cdots \int_0^{x-\sum_{i=1}^{n-2} u_i} \left(\prod_{i=1}^{n-1}(\theta + u_i^2)\right)\left(\theta + x - \sum_{i=1}^{n-1} u_i^2\right) du_{n-1} \ldots du_2 du_1 = \mathrm{D}_{n-1}.$$

So, $f_{S_n}(x)$ with Pranav distribution is equals to

$$f_{S_n}(x) = \left(\frac{\theta^4}{\theta^4 + 6}\right)^n e^{-\theta x} \sum_{r=0}^{n}\left(\binom{n}{r} 6^r \theta^{n-r} \frac{x^{n+3r-1}}{(n+3r-1)!}\right). \tag{7}$$

**2.5. Rani distribution**

Replacement the Rani distribution in Eq. (3) and for $x > 0$, the functions of $f_{S_n}(x)$, for $n=2$ and 3 are given by

$$f_{S_2}(x) = \left(\frac{\theta^5}{\theta^5 + 24}\right)^2 e^{-\theta x} \int_0^x (\theta + u_1^4)(\theta + (x-u_1)^4) du_1 = \left(\frac{\theta^5}{\theta^5 + 24}\right)^2 e^{-\theta x} \mathrm{E}_1,$$

$$f_{S_3}(x) = \left(\frac{\theta^5}{\theta^5 + 24}\right)^3 e^{-\theta x} \int_0^x \int_0^{x-u_1} (\theta + u_1^4)(\theta + u_2^4)(\theta + (x-u_1-u_2)^4) du_2 du_1$$

$$= \left(\frac{\theta^5}{\theta^5 + 24}\right)^3 e^{-\theta x} \mathrm{E}_2.$$

where

$$\mathrm{E}_1 = \theta^2 x + \frac{2\theta x^5}{5} + \frac{x^9}{630},$$

$$\mathrm{E}_2 = \frac{\theta^3 x^2}{2} + \frac{\theta^2 x^6}{10} + \frac{\theta x^{10}}{2100} + \frac{x^{14}}{6306300}.$$

Subsequently

$$\mathrm{E}_1 = \binom{2}{0} 24^0 \theta^{2-0} \frac{x^1}{1!} + \binom{2}{1} 24^1 \theta^{2-1} \frac{x^5}{5!} + \binom{2}{2} 24^2 \theta^{2-2} \frac{x^9}{9!},$$

$$\mathrm{E}_2 = \binom{3}{0} 24^0 \theta^{3-0} \frac{x^2}{2!} + \binom{3}{1} 24^1 \theta^{3-1} \frac{x^6}{6!} + \binom{3}{2} 24^2 \theta^{3-2} \frac{x^{10}}{10!} + \binom{3}{3} 24^3 \theta^{3-3} \frac{x^{14}}{141!}.$$

According to the above results, it can be concluded that $\mathrm{E}_{n-1}$ is equals to

$$\mathrm{E}_{n-1} = \sum_{r=0}^{n} \binom{n}{r} 24^r \theta^{n-r} \frac{x^{n+4r-1}}{(n+4r-1)!}.$$

Thus

$$\int_0^x \int_0^{x-u_1} \cdots \int_0^{x-\sum_{i=1}^{n-2} u_i} \left(\prod_{i=1}^{n-1}(\theta + u_i^2)\right)\left(\theta + x - \sum_{i=1}^{n-1} u_i^2\right) du_{n-1} \ldots du_2 du_1 = \mathrm{E}_{n-1}.$$

So, $f_{S_n}(x)$ with Rani distribution is equals to

$$f_{S_n}(x) = \left(\frac{\theta^5}{\theta^5 + 24}\right)^n e^{-\theta x} \sum_{r=0}^{n}\left(\binom{n}{r} 24^r \theta^{n-r} \frac{x^{n+4r-1}}{(n+4r-1)!}\right). \tag{8}$$

### 2.6. Ram Awadh distribution

Replacement the Rani distribution in Eq. (3) and for $x > 0$, the functions of $f_{S_n}(x)$, for $n=2$ and 3 are given by

$$f_{S_2}(x) = \left(\frac{\theta^6}{\theta^6 + 120}\right)^2 e^{-\theta x} \int_0^x (\theta + u_1^5)(\theta + (x - u_1)^5) du_1 = \left(\frac{\theta^6}{\theta^6 + 120}\right)^2 e^{-\theta x} F_1,$$

$$f_{S_3}(x) = \left(\frac{\theta^6}{\theta^6 + 120}\right)^3 e^{-\theta x} \int_0^x \int_0^{x-u_1} (\theta + u_1^5)(\theta + u_2^5)(\theta + (x - u_1 - u_2)^5) du_2 du_1$$

$$= \left(\frac{\theta^6}{\theta^6 + 120}\right)^3 e^{-\theta x} F_2.$$

where

$$F_1 = \theta^2 x + \frac{\theta x^6}{3} + \frac{x^{11}}{2772},$$

$$F_2 = \frac{\theta^3 x^2}{2} + \frac{\theta^2 x^7}{14} + \frac{\theta x^{12}}{11088} + \frac{x^{17}}{205837632}.$$

Subsequently

$$F_1 = \binom{2}{0} 120^0 \theta^{2-0} \frac{x^1}{1!} + \binom{2}{1} 120^1 \theta^{2-1} \frac{x^6}{6!} + \binom{2}{2} 120^2 \theta^{2-2} \frac{x^{11}}{11!},$$

$$F_2 = \binom{3}{0} 120^0 \theta^{3-0} \frac{x^2}{2!} + \binom{3}{1} 120^1 \theta^{3-1} \frac{x^7}{7!} + \binom{3}{2} 120^2 \theta^{3-2} \frac{x^{12}}{12!} + \binom{3}{3} 120^3 \theta^{3-3} \frac{x^{17}}{17!}.$$

According to the above results, it can be concluded that $F_{n-1}$ is equals to

$$F_{n-1} = \sum_{r=0}^n \binom{n}{r} 120^r \theta^{n-r} \frac{x^{n+5r-1}}{(n+5r-1)!}.$$

Thus

$$\int_0^x \int_0^{x-u_1} \cdots \int_0^{x-\sum_{i=1}^{n-2} u_i} \left(\prod_{i=1}^{n-1}(\theta + u_i^2)\right)\left(\theta + x - \sum_{i=1}^{n-1} u_i^2\right) du_{n-1} \ldots du_2 du_1 = F_{n-1}.$$

So, $f_{S_n}(x)$ with Rani distribution is equals to

$$f_{S_n}(x) = \left(\frac{\theta^6}{\theta^6 + 120}\right)^n e^{-\theta x} \sum_{r=0}^n \left(\binom{n}{r} 120^r \theta^{n-r} \frac{x^{n+5r-1}}{(n+5r-1)!}\right). \tag{9}$$

### 3. The *m*th Moments

Another important measure for statistical distributions is the *m*th moment. The calculation of this factor is very important for any statistical distribution. The *m*th moment for Eqs. (4)-(9), can be expressed as

$$E(X^m) = \int_0^\infty x^m f_{S_n}(x) \, dx. \tag{10}$$

By special values for amounts of *m* can be defined the descriptive statistical indicators such as, mean, variance, coefficient of variation, skewness and kurtosis. Table 3 presents the closed form for the *m*th moments of PDFs Table 2 as follows.

**Table 3.** The *m*th moments for PDFs of Table 2

| Distributions | *m*th moments |
|---|---|
| Shanker | $\frac{m!}{\theta^m} \left( \frac{\theta^2}{\theta^2+1} \right)^n \sum_{r=0}^{n} \binom{n}{r} \binom{n+m+r-1}{n+r-1} \left( \frac{1}{\theta^2} \right)^r$ |
| Akash | $\frac{m!}{\theta^m} \left( \frac{\theta^2}{\theta^2+2} \right)^n \sum_{r=0}^{n} \binom{n}{r} \binom{n+m+2r-1}{n+2r-1} \left( \frac{2}{\theta^2} \right)^r$ |
| Ishita | $\frac{m!}{\theta^m} \left( \frac{\theta^3}{\theta^3+2} \right)^n \sum_{r=0}^{n} \binom{n}{r} \binom{n+m+2r-1}{n+2r-1} \left( \frac{2}{\theta^3} \right)^r$ |
| Pranav | $\frac{m!}{\theta^m} \left( \frac{\theta^4}{\theta^4+6} \right)^n \sum_{r=0}^{n} \binom{n}{r} \binom{n+m+3r-1}{n+3r-1} \left( \frac{6}{\theta^4} \right)^r$ |
| Rani | $\frac{m!}{\theta^m} \left( \frac{\theta^5}{\theta^5+24} \right)^n \sum_{r=0}^{n} \binom{n}{r} \binom{n+m+4r-1}{n+4r-1} \left( \frac{24}{\theta^5} \right)^r$ |
| Ram Awadh | $\frac{m!}{\theta^m} \left( \frac{\theta^6}{\theta^6+120} \right)^n \sum_{r=0}^{n} \binom{n}{r} \binom{n+m+5r-1}{n+5r-1} \left( \frac{120}{\theta^6} \right)^r$ |

## 4. Application

As mentioned earlier, the whole lifetime for a cold standby spare systems 1 out of *n* is equal to $S_n = X_1 + X_2 + \cdots + X_n = \sum_{i=1}^{n} X_i$. In addition, it has been proved that the Lindley distribution in many aspects, such as the systems reliability is more flexible than the exponential distribution. Therefore, with the knowledge of $f_{S_n}(x)$ by Lindley distribution (the proof of it is given in Appendix 1.), can be obtained a general formula for the reliability of the cold standby spare system under failure time the Lindley.

The calculation of the reliability for a cold standby redundant system at time *t* ($R_S(t)$), with a prime active unit and (*n*-1)-spare unit according to (Coit (2001)), is equal to

$$R_S(t) = P(S_n > t) = P(X_1 + \cdots + X_n > t) = r(t) + \sum_{i=1}^{n-1} \int_0^t f_{S_i}(u) r(t-u) \, du, \tag{11}$$

where $r(t)$ is reliability of a component at time *t*, and $f_{S_i}(t)$ is the pdf for sums of *i*-unit failure time.

Nevertheless, after some mathematical calculations, the system reliability when the failure behavior of its components follows the Lindley distribution is determined as

$$R(t)_{\text{Lindley}} = e^{-\theta t} \left( \sum_{i=0}^{n-1} \sum_{j=0}^{i} \left(\frac{\theta^2}{1+\theta}\right)^i \binom{i}{j} \frac{t^{i+j}}{(i+j)!} + \frac{\theta}{1+\theta} \sum_{i=0}^{n-1} \sum_{j=0}^{i} \left(\frac{\theta^2}{1+\theta}\right)^i \binom{i}{j} \frac{t^{i+j+1}}{(i+j+1)!} \right). \quad (12)$$

The system mean time to failure (MTTF) under the Lindley distribution is also equal

$$\text{MTTF}_{\text{Lindley}} = \int_0^\infty R(t)_{\text{Lindley}} dt$$

$$= \left( \sum_{i=0}^{n-1} \sum_{j=0}^{i} \frac{\binom{i}{j}\left(\frac{\theta^2}{1+\theta}\right)^i}{(i+j)!} \int_0^\infty t^{i+j} e^{-\theta t} dt + \frac{\theta}{1+\theta} \sum_{i=0}^{n-1} \sum_{j=0}^{i} \frac{\binom{i}{j}\left(\frac{\theta^2}{1+\theta}\right)^i}{(i+j+1)!} \int_0^\infty t^{i+j+1} e^{-\theta t} dt \right). \quad (13)$$

Simplicity can be prove that

$$\int_0^\infty t^{i+j} e^{-\theta t} dt = \frac{(i+j)!}{\theta^{i+j+1}}, \int_0^\infty t^{i+j+1} e^{-\theta t} dt = \frac{(i+j+1)!}{\theta^{i+j+2}}.$$

So, Eq. (13) can be expressed as

$$\text{MTTF}_{\text{Lindley}} = \int_0^\infty R(t)_{\text{Lindley}} dt = \frac{2+\theta}{\theta(1+\theta)} \sum_{i=0}^{n-1} \sum_{j=0}^{i} \binom{i}{j} \left(\frac{\theta}{1+\theta}\right)^i \left(\frac{1}{\theta}\right)^j. \quad (14)$$

By using the binomial expansion, we know very well that $\sum_{j=0}^{i} \binom{i}{j} \left(\frac{1}{\theta}\right)^j = \left(1 + \frac{1}{\theta}\right)^i$. Therefore, Eq. (14) reduces to

$$\text{MTTF}_{\text{Lindley}} = \frac{n}{\theta}\left(\frac{2+\theta}{1+\theta}\right). \quad (15)$$

Furthermore, the general equation for the system reliability and MTTF with exponential distribution are calculated as (Billinton & Allan (1992))

$$R(t)_{\text{Exponential}} = e^{-\theta t} \sum_{i=0}^{n-1} \frac{(\theta t)^i}{i!}, \quad (16)$$

$$\text{MTTF}_{\text{Exponential}} = \frac{n}{\theta}. \quad (17)$$

Now, by setting $\theta = 0.1, 0.5, 1$ and $3$, and mission time 100 hours, we are interested to compare Eqs. (12) and (16) by 5-component in the system (i.e., $n=5$).

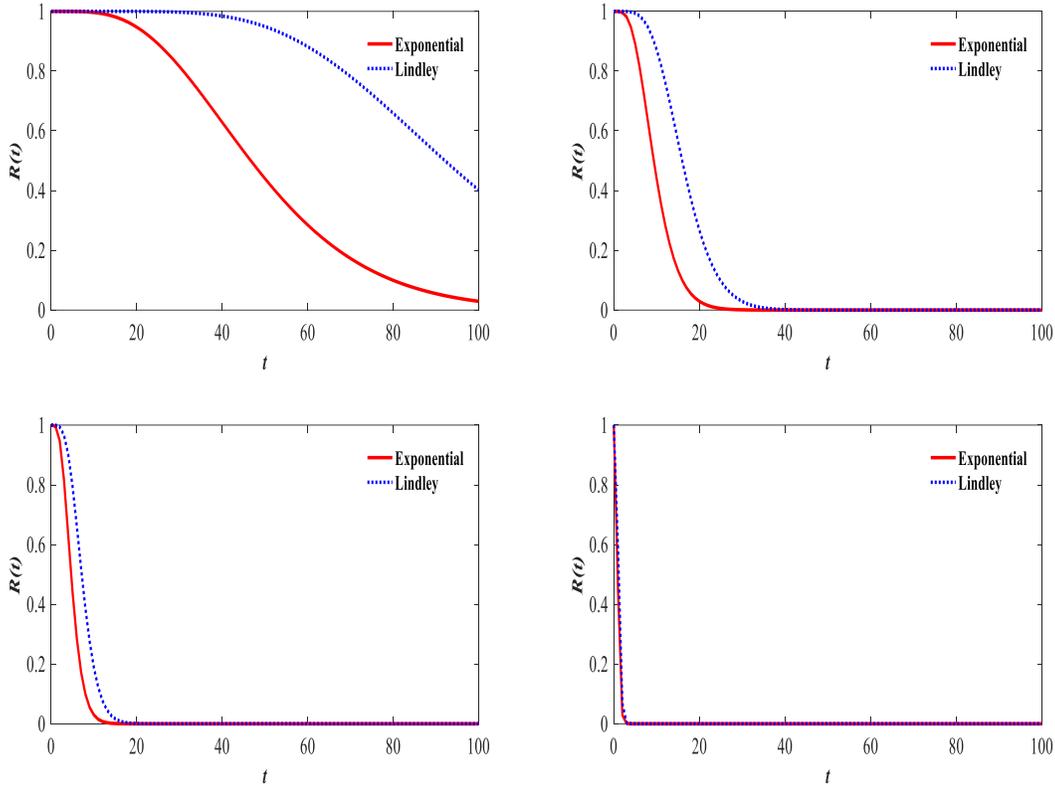

**Figure 1.** Systems reliability

In Figure 1 as expected, the system reliability under the Lindley distribution is always higher than the system reliability under the exponential distribution. Moreover, according to the Eqs. (15) and (17), for different values of $\theta$ that already mentioned, The MTTFs are given by Table 4.

**Table 4.** The comparison of MTTFs

| $\theta$ | 0.1 | 0.5 | 1 | 3 |
|---|---|---|---|---|
| $MTTF_{Lindley}$ | 95.45 | 16.67 | 7.5 | 2.08 |
| $MTTF_{Exponential}$ | 40 | 10 | 5 | 1.67 |

The results of the Table 4 just like Figure 1, show that the average time for the system during operation always performs better under the Lindley distribution.

5. **Conclusion**

In this paper, the probability density functions of the sum of *n* independent and identically distributed random variables with Shanker, Akash, Ishita, Pranav, Rani, and Ram Awadh have been explicitly assessed. Their determination is based on the change-of-variables technique.

Moreover, the *m*th moments for them were also accurately calculated. Finally, as an application, a general formula was calculated for a 1 out of *n* cold standby spare system reliability with component failure under the Lindley distribution.

**Appendix 1.**

Using Eq. (3), for amounts of *n*=2 and 3 functions of $f_{S_n}(x)$, for $x > 0$ are given as

$$f_{S_2}(x) = \left(\frac{\theta^2}{1+\theta}\right)^2 e^{-\theta x} \int_0^x (1+u_1)(1+x-u_1)du_1 = \left(\frac{\theta^2}{1+\theta}\right)^2 e^{-\theta x} L_1,$$

$$f_{S_3}(x) = \left(\frac{\theta^2}{1+\theta}\right)^3 e^{-\theta x} \int_0^x \int_0^{x-u_1} (1+u_1)(1+u_2)(1+x-u_1-u_2)du_2 du_1 = \left(\frac{\theta^2}{1+\theta}\right)^3 e^{-\theta x} L_2.$$

where

$$L_1 = x + x^2 + \frac{x^3}{6},$$

$$L_2 = \frac{x^2}{2} + \frac{x^3}{2} + \frac{x^4}{8} + \frac{x^5}{120}.$$

Subsequently, after some algebra, we have

$$L_1 = \frac{x^1}{1!} + 2\frac{x^2}{2!} + \frac{x^3}{3!} = \binom{2}{0}\frac{x^1}{1!} + \binom{2}{1}\frac{x^2}{2!} + \binom{2}{2}\frac{x^3}{3!},$$

$$L_2 = \frac{x^2}{2!} + 3\frac{x^3}{3!} + 3\frac{x^4}{4!} + \frac{x^5}{5!} = \binom{3}{0}\frac{x^2}{2!} + \binom{3}{1}\frac{x^3}{3!} + \binom{3}{2}\frac{x^4}{4!} + \binom{3}{3}\frac{x^5}{51!}.$$

As can be seen, the coefficients in the expressions $L_1$ and $L_2$ make the Pascal's triangle. According to the above results, it can be concluded that $L_{n-1}$ is equal to

$$L_{n-1} = \sum_{r=0}^{n} \binom{n}{r} \frac{x^{n+r-1}}{(n+r-1)!}.$$

Thus

$$\int_0^x \int_0^{x-u_1} \cdots \int_0^{x-\sum_{i=1}^{n-2} u_i} \left(\prod_{i=1}^{n-1}(1+u_i)\right)\left(1+x-\sum_{i=1}^{n-1} u_i\right) du_{n-1} \cdots du_2 du_1 = L_{n-1}.$$

So, $f_{S_n}(x)$ with Lindley distribution is equals to

$$f_{S_n}(x) = \left(\frac{\theta^2}{\theta+1}\right)^n e^{-\theta x} \sum_{r=0}^{n} \binom{n}{r} \frac{x^{n+r-1}}{(n+r-1)!}. \tag{18}$$

where Eq. (18) is exactly the same formula that was presented by Zakerzadeh & Dolati (2009), Al-Mutairi et al. (2013) and Hassan (2014).